# SEMISTABLE ABELIAN VARIETIES WITH SMALL DIVISION FIELDS

ARMAND BRUMER AND KENNETH KRAMER

1. INTRODUCTION

The conjecture of Shimura-Taniyama-Weil, now proved through the work of Wiles and disciples, is only part of the Langlands program. Based on a comparison of the local factors ([And], [Ser1]), it also predicts that the L-series of an abelian surface defined over $\mathbb{Q}$ should be the L-series of a Hecke eigen cusp form of weight 2 on a suitable group commensurable with $\mathrm{Sp}_4(\mathbb{Z})$. The only decisive examples are related to lifts of automorphic representations of proper subgroups of $\mathrm{Sp}_4$, for example the beautiful work of Yoshida ([Yos], [BS1]).

We do not know how to construct "non-trivial" weight 2 forms for groups such as the paramodular group or the Hecke subgroup $\Gamma_0(p)$ in genus 2. In fact, it seems that not a single example is known which is not a lift. Not even a conjectural dimension formula for weight 2 cusp forms has yet been proposed ([Ibu],[Has]).

It seems interesting to consider existence and non-existence on the arithmetic side. If the special fiber $\mathcal{A}_p$ of the Néron model of the abelian variety $A$ has no unipotent part, we say that $A$ has semistable reduction at $p$. In that case, the exponent $\mathfrak{f}_p(A)$ of $p$ in the conductor of $A$ is the toroidal dimension of $\mathcal{A}_p$. For example, the modular variety $J_0(p)$ has conductor $p^d$, where $d$ is its dimension. To ensure that the putative automorphic form not be a lift, we could study surfaces of prime conductor. This guarantees that we are not dealing, for instance, with a surface isogenous to a product of elliptic curves, the Weil restriction of scalars of an elliptic curve defined over a quadratic number field or a surface with non-trivial endomorphisms. Serre's conjecture [Ser2] implies that such a surface should be $\mathbb{Q}$-isogenous to the Jacobian of a curve also defined over $\mathbb{Q}$. This led us to try to understand curves $C$ of genus 2 whose Jacobian variety $J = J(C)$ has prime conductor $p$. A minimal model for $C$ has bad reduction at $p$, but also may have primes of bad reduction where $J$ has good reduction. We call such primes, at which $J$ reduces to the product of 2 elliptic curves, primes of *mild* reduction. Because of mild primes, we could not immediately rule out the possibility that all six Weierstrass points of $C$ be rational, that is $\mathbb{Q}(J[2]) = \mathbb{Q}$. In the present investigation, we find that $\mathrm{Gal}(\mathbb{Q}(J[2])/\mathbb{Q})$ cannot even be nilpotent, as we shall explain below.

One expects that $\mathrm{Gal}(\mathbb{Q}(A[\ell])/\mathbb{Q})$ should be, in general, as large as is compatible with the endomorphisms of $A$ and the Weil pairing. When $\mathrm{End}\, A = \mathbb{Z}$, this has been proved for $\ell$ sufficiently large by Serre [Ser3], in all odd dimensions as well as in dimensions 2 and 6. Note that $\mathbb{Q}(A[\ell])$ always contains the cyclotomic field $\mathbb{Q}(\boldsymbol{\mu}_\ell)$,

*Date*: June 20, 2002.

Work of the second author was supported in part by a grant from the City University of New York PSC-CUNY Research Award Program.

whether or not $A$ is principally polarized. As examples of the results obtained here, we have:

**Proposition 3.5.** *Let $A/\mathbb{Q}$ be a semistable abelian variety with good reduction at a prime $\ell \leq 19$ and set $H = \mathrm{Gal}(\mathbb{Q}(A[\ell])/\mathbb{Q}(\boldsymbol{\mu}_\ell))$. Then $\ell$ annihilates the maximal abelian quotient $H/H'$. So $H$ is an $\ell$-group if it is nilpotent.*

**Theorem 1.1.** *Suppose that $A/\mathbb{Q}$ is semistable, with bad reduction only at $p$, and that $H$ is nilpotent. If $\ell = 2$ or $3$, then $\mathfrak{f}_p(A) = \dim A$, so the reduction at $p$ is totally toroidal. Moreover $p \equiv 1 \bmod 8$ if $\ell = 2$ and $p \equiv 1 \bmod 3$ if $\ell = 3$.*

Consider elliptic curves of prime conductor $p$ with some rational $\ell$-torsion. For $\ell = 2$, Neumann [Ne] and Setzer [Sz] found that such elliptic curves exist if and only if $p = n^2 + 64$ or $p = 17$. For odd $\ell$, Miyawaki [Mi] showed that the well-known examples with $(\ell, p) = (3, 19)$, $(3,37)$ or $(5,11)$ are the only ones. Under the stronger assumption that $A$ has prime conductor, Theorem 1.1 implies:

**Corollary 1.2.** *Let $A/\mathbb{Q}$ be an abelian variety of prime conductor $p$ such that $H$ is nilpotent. If $\ell = 2$, then $A$ is a Neumann-Setzer elliptic curve. If $\ell = 3$, then $p = 19$ or $37$ and $A$ is a Miyawaki elliptic curve.*

Let $\mathcal{C}$ be the $\mathbb{Q}$-isogeny class of an abelian variety $A$ with semistable bad reduction at $p$. We denote by $\Phi_A$ the group of connected components of the special fiber of its Néron model at $p$. We say that $A$ is $\ell$-maximal if $\mathrm{ord}_\ell |\Phi_A(\overline{\mathbb{F}}_p)|$ is maximal among varieties in $\mathcal{C}$. A major role will be played by the collection

$$(1.3) \qquad \mathcal{V}_\ell(\mathcal{C}) = \{A \in \mathcal{C} \mid \mathrm{ord}_\ell |\Phi_A(\overline{\mathbb{F}}_p)| \text{ is maximal}\}$$

of $\ell$-maximal varieties. This notion depends on the choice of $p$, which will be clear from the context.

We show that $\mathbb{Q}(A[\ell]) = \mathbb{Q}(\boldsymbol{\mu}_{2\ell})$ for all $A \in \mathcal{V}_\ell(\mathcal{C})$ if $\mathcal{C}$ is the isogeny class of a product of elliptic curves of conductor $p$ with some rational $\ell$-torsion, cf. Proposition 6.7. The following is a converse when $\ell = 2$ or $3$.

**Theorem 1.4.** *Let $\ell = 2$ or $3$ and let $\mathcal{C}$ be a $\mathbb{Q}$-isogeny class of semistable abelian varieties with bad reduction only at $p$. If $\mathbb{Q}(A[\ell]) \subseteq \mathbb{Q}(\boldsymbol{\mu}_{2\ell})$ for all $A \in \mathcal{V}_\ell(\mathcal{C})$, then $\mathcal{C}$ is the isogeny class of a product of Neumann-Setzer ($\ell = 2$) or Miyawaki ($\ell = 3$) curves of conductor $p$.*

A brief outline of the paper follows. With the help of class field theory and the conductor bounds of Fontaine, we obtain in section 3 the information we require about the group $H = \mathrm{Gal}(\mathbb{Q}(A[\ell])/\mathbb{Q}(\boldsymbol{\mu}_\ell))$.

Denote by $\mathbb{T}_\ell(A)$ the Tate module of $A$ and by $\mathcal{D}_v$ the decomposition group in $G_\mathbb{Q}$ of a place $v$ over $p$. In section 4, we recall the definition of the $\mathcal{D}_v$-submodules of $\mathbb{T}_\ell(A)$ introduced by Grothendieck [Gro]. Now suppose $A$ is semistable with bad reduction only at $p$ and that $H$ is an $\ell$-group. Let $a$ denote the dimension of the abelian variety in the special fiber of $A$ at $p$. We use the Grothendieck modules together with $\ell$-maximality and a strange lifting lemma to construct a pure $\mathbb{Z}_\ell$-submodule $W_A$ of $\mathbb{T}_\ell(A)$ of rank $2a$. Although $W_A$ is canonically only a $\mathcal{D}_v$-module, we obtain the rather surprising result that it is a $G_\mathbb{Q}$-module under suitable hypotheses. Then, in section 5, we use the formal group associated to $A$ at $\ell$ and Raynaud's theory of $\mathbb{F}$-vector space schemes to prove that $W_A = 0$, so that $a = 0$ and $A$ is totally toroidal at $p$. The proof of Theorem 1.1 may be found at the end of this section.

A generalization of Theorem 1.4 is proved in section 6, using the theorem of Faltings on Tate's conjecture. We then give two examples for which $H$ is an $\ell$-group, namely, $J_0(31)$ with $\ell = 2$ and $J_0(41)$ with $\ell = 5$. Since they clearly are not isogenous to a product of elliptic curves, any attempt to weaken the hypotheses of Theorem 1.4 will require some care.

Fontaine [Fo] proved that there are no abelian varieties of positive dimension with everywhere good reduction over $\mathbb{Q}$. We showed, in [BK2], that no semistable abelian variety exists with bad reduction exactly at one prime $p \leq 7$ and Schoof [Scho] improved this result by different methods.

In our next paper, using other techniques, we plan to study the existence of abelian surfaces $A/\mathbb{Q}$ of conductor $p$ in the spirit of [BK1]. When $A$ is principally polarized, $G = \text{Gal}(\mathbb{Q}(A[2])/\mathbb{Q})$ is a subgroup of $\text{GSp}_4(\mathbb{F}_2) \simeq \mathcal{S}_6$, the symmetric group on 6 letters. As an application of Theorem 1.1, we will show that $G$ can only be $\mathcal{S}_n$ for $3 \leq n \leq 6$, $\mathcal{S}_n \times \mathbb{Z}/2$ for $n = 3, 4$, or the wreath product $\mathcal{S}_3 \wr \mathbb{Z}/2$, with each of these possibilities occurring infinitely often under the Schinzel hypothesis [SS].

Our criteria will show that $p = 277$ is the smallest prime for which there is a curve of genus 2 whose Jacobian has prime conductor. One such curve is given by
$$y^2 + y = x^5 - 2x^3 + 2x^2 - x.$$
We thus propose, as a challenge to our automorphic friends, a verification that all modular forms of weight 2 for smaller prime level are "lifts", but that at least one interesting eigenform exists for 277.

This paper concentrates on one of the results presented by the first author at this Conference. He is grateful to K. Hashimoto and K. Miyake for their kind invitation and for their support both moral and financial. Other results, more influenced by the work of Hashimoto and Ibukiyama, will be published later. Thanks to H. Nakamura for helping run the Conference so pleasantly and seamlessly.

## 2. Acknowledgments

The first author wishes to express his indebtedness for the generous hospitality, patience, friendship and stimulation afforded him by Professor Hashimoto and his students at Waseda University. He also thanks Professors K. Miyake of Tokyo Metropolitan University, T. Ibukiyama of Osaka University and H. Saito of Kyoto University, their institutions and colleagues. Together, they made it possible to spend a wonderful month in Japan, beyond the week of this Conference.

## 3. Controlling the $\ell$-division field

Let $S$ be a set of rational primes not containing the prime $\ell$ nor the archimedean prime $\infty$. We shall say that a finite extension $L/\mathbb{Q}$ is $(\ell, S)$-controlled if it has the following properties:

(L1) $L$ is Galois over $\mathbb{Q}$ and contains $\boldsymbol{\mu}_\ell$;
(L2) $L/\mathbb{Q}$ is unramified outside each place over $S \cup \{\ell, \infty\}$;
(L3) the ramification degree for each place of $L$ over $S$ divides $\ell$;
(L4) the higher ramification groups $\mathcal{D}_\lambda^{(u)}$ are trivial for $u > 1/(\ell - 1)$, where $\mathcal{D}_\lambda$ is the decomposition group at a place $\lambda$ over $\ell$ in $L$.

If $S = \{p\}$ consists of one prime, we say $L/\mathbb{Q}$ is $(\ell, p)$-controlled. Thanks to Grothendieck [Gro, §2.5] and Fontaine [Fo, Thm. A], the $\ell$-division field of a semistable abelian variety $A/\mathbb{Q}$ with good reduction at $\ell$ is $(\ell, S)$-controlled when $S$ contains the primes of bad reduction. See also [BK2, §3].

Since the upper ramification numbering behaves well with respect to quotients, it is clear that a subfield of $L$ containing $\boldsymbol{\mu}_\ell$ and Galois over $\mathbb{Q}$ also is $(\ell, S)$-controlled. Properties (L1) and (L2) certainly are preserved under composition of $(\ell, S)$-controlled extensions. By uniqueness of the tamely ramified extension of degree $\ell$ over the unramified closure of $\mathbb{Q}_p$, property (L3) is preserved. Property (L4) is inherited by the compositum of finitely many Galois extensions that enjoy it, again using the fact that the upper ramification is preserved by quotients.

It is helpful to highlight the impact of (L4) on the following local situation, a variant of which was treated in [BK2, Lem. 6].

**Lemma 3.1.** *Let $E \supseteq F \supseteq \mathbb{Q}_\ell(\boldsymbol{\mu}_\ell) \supseteq \mathbb{Q}_\ell$ be a tower of local fields such that $E/\mathbb{Q}_\ell$ is Galois and $E/\mathbb{Q}_\ell(\boldsymbol{\mu}_\ell)$ is an abelian $\ell$-extension. The higher ramification groups $\operatorname{Gal}(E/\mathbb{Q}_\ell)^{(u)}$ are trivial for all $u > 1/(\ell - 1)$ if and only if each of the abelian conductor exponents $\mathfrak{f}(E/F)$ and $\mathfrak{f}(F/\mathbb{Q}_\ell(\boldsymbol{\mu}_\ell))$ is at most 2.*

*Proof.* First consider a general Galois extension of local fields $E/K$ and put $H = \operatorname{Gal}(E/K)$. Assume only that $F$ is the fixed field of a normal subgroup $N \triangleleft H$ and let $\bar{H} = H/N$. Write $\varphi_{E/K}$ for the Herbrand function [Se, IV, §3] from the lower to upper ramification numbering of $H$. Let us verify that $H_x = 1$ if and only if both $N_x$ and $\bar{H}_{x'}$ are trivial, where $x' = \varphi_{E/F}(x)$. If $H_x = 1$, then $N_x = H_x \cap N = 1$ and by Herbrand's theorem [Se, IV, §3, Lem. 5], we also have $\bar{H}_{x'} = H_x N/N = 1$. Conversely, if $\bar{H}_{x'} = 1$, then $H_x \subseteq N$. It follows that $H_x = N_x$, so if $N_x$ also is trivial, we may conclude that $H_x = 1$. This observation will be used below with $K = \mathbb{Q}_\ell(\boldsymbol{\mu}_\ell)$.

Now assume all the given hypotheses and let $D = \operatorname{Gal}(E/\mathbb{Q}_\ell)$. Since the tame ramification degree of the extension $E/\mathbb{Q}_\ell$ is $\ell - 1$, we have $\varphi_{E/\mathbb{Q}_\ell}(1) = |D_1|/|D_0| = 1/(\ell - 1)$. Furthermore, $D_x = H_x$ for all $x \geq 1$ because $H$ is the $\ell$-Sylow subgroup of $D$. Then we have the following chain of equivalences:

$$\begin{aligned} D^{(u)} = 1 \text{ for all } u > 1/(\ell-1) &\Leftrightarrow D_{1+\epsilon} = 1 \text{ for all } \epsilon > 0 \\ &\Leftrightarrow H_{1+\epsilon} = 1 \text{ for all } \epsilon > 0 \\ &\Leftrightarrow N_{1+\epsilon} = 1 \text{ and } \bar{H}_{1+\epsilon} = 1 \text{ for all } \epsilon > 0, \end{aligned}$$

using our opening observation and the fact that $\varphi_{E/F}(1) = 1$ for the last equivalence.

The conductor exponent of an abelian extension $E/K$ is given by $\mathfrak{f}(E/K) = \varphi_{E/K}(c) + 1$, where $c$ is the largest integer such that $H_c \neq 1$. When $H$ is an $\ell$-group, we have $\varphi_{E/K}(1) = 1$. For an abelian $\ell$-extension $E/K$ we therefore find that $\mathfrak{f}(E/K) \leq 2$ if and only if $H_{1+\epsilon} = 1$ for all $\epsilon > 0$. Applying this fact to the abelian extensions $E/F$ and $F/\mathbb{Q}_\ell(\boldsymbol{\mu}_\ell)$, we may conclude that $D^{(u)} = 1$ for all $u > 1/(\ell-1)$ if and only if both $\mathfrak{f}(E/F) \leq 2$ and $\mathfrak{f}(F/\mathbb{Q}_\ell(\boldsymbol{\mu}_\ell)) \leq 2$. □

We use below the fact, proved in [MM], that $\mathbb{Q}(\boldsymbol{\mu}_\ell)$ has class number one exactly when $\ell \leq 19$.

**Lemma 3.2.** *Suppose $M$ is $(\ell, S)$-controlled and unramified outside $\ell$, with $\ell \leq 19$. If $M$ is abelian over $F = \mathbb{Q}(\boldsymbol{\mu}_{2\ell})$, then $M = F$.*

*Proof.* To analyze the extension $M/F$, we introduce the usual notation of class field theory. Let $\mathbb{A}_F^\times$ be the idele group of $F$ and write $\mathbf{F}^\times$ for the image of $F^\times$ on the diagonal. Let $F_v$ be the completion of $F$ at a place $v$ and put $U_v$ for the connected component of the identity in $F_v^\times$. If $v$ is non-archimedean, set $U_v^{(m)} = \{x \in U_v \mid v(x-1) \geq m\}$.

Observe that there is a unique prime $\lambda$ over $\ell$ in $F$ and let $\mathfrak{l}$ be a prime of $M$ over $\lambda$. By (L4), we have $\mathrm{Gal}(M_\mathfrak{l}/\mathbb{Q}_\ell)^{(u)} = 1$ for all $u > 1/(\ell-1)$. For brevity, write $\mathfrak{f}_\lambda$ for the abelian conductor exponent of the local extension $M_\mathfrak{l}/F_\lambda$. Since $M/F$ is unramified outside $\ell$, there is a surjection $\mathbb{A}_F^\times/\mathcal{N} \twoheadrightarrow \mathrm{Gal}(M/F)$, where $\mathcal{N}$ is the norm subgroup

$$\mathcal{N} = (U_\lambda^{(\mathfrak{f}_\lambda)} \times \prod_{v \neq \lambda} U_v) \cdot \mathbf{F}^\times.$$

Let $V_F$ be the image of the global units of $F$ in $\Gamma = U_\lambda/U_\lambda^{(\mathfrak{f}_\lambda)}$. In our case $U_v \simeq \mathbb{C}^\times$ for archimedean places, so we have the exact sequence

(3.3) $$0 \to \Gamma/V_F \to \mathbb{A}_F^\times/\mathcal{N} \to \mathbb{A}_F^\times/((\prod_v U_v) \cdot \mathbf{F}^\times) \to 0,$$

in which the last quotient is isomorphic to the ideal class group of $F$ and is trivial by assumption.

The group $U_\lambda/U_\lambda^{(1)}$ is generated by the images of global cyclotomic units when $\ell$ is odd and is trivial when $\ell = 2$. Since $U_\lambda^{(1)}/U_\lambda^{(\mathfrak{f}_\lambda)}$ is an $\ell$-group, we find that $\Gamma/V_F$ is an $\ell$-group and so $\mathrm{Gal}(M/F)$ is an $\ell$-group by exact sequence (3.3). Then $\mathfrak{f}_\lambda \leq 2$ by Lemma 3.1. The group $U_\lambda^{(1)}/U_\lambda^{(2)}$ is generated by the image of a global primitive $\ell^{\mathrm{th}}$ root of unity when $\ell$ is odd or by the image of $i$ when $\ell = 2$. Thus $\Gamma/V_F$ and *a fortiori* $\mathrm{Gal}(M/F)$ is trivial. □

**Remark.** Under the hypothesis that $\mathrm{Gal}(M/F)$ is an abelian $\ell$-group, similar reasoning shows that $M = F$ when $\ell$ is a regular prime.

**Lemma 3.4.** *Suppose that $\ell \leq 19$ and that $L$ is an $(\ell, S)$-controlled extension containing $F = \mathbb{Q}(\boldsymbol{\mu}_{2\ell})$. Let $E$ be the maximal subfield of $L$ abelian over $F$. Then $\mathrm{Gal}(E/F)$ is annihilated by $\ell$ and $\dim_{\mathbb{F}_\ell} \mathrm{Gal}(E/F) \leq s$, where $s$ is the number of primes over $S$ in $F$.*

*Proof.* For each prime $\mathfrak{p}_j$ over $p$ in $F$, choose a prime $\mathfrak{P}_j$ over $\mathfrak{p}_j$ in $L$. Since $p$ does not ramify in $F/\mathbb{Q}$, the inertia group $\mathcal{I}_j = \mathcal{I}_{\mathfrak{P}_j}$ is contained in $H = \mathrm{Gal}(L/F)$. Let $N$ be the subgroup of $H$ generated by $\mathcal{I}_j$ for $j = 1, \ldots, s$ and by the commutator subgroup $H'$. The previous Lemma shows that $N = H$. By (L3), each $\mathcal{I}_j$ has order dividing $\ell$, so $\mathrm{Gal}(E/F)$ is annihilated by $\ell$ and its $\mathbb{F}_\ell$-dimension is at most $s$. □

**Remark.** The bound on $\dim_{\mathbb{F}_\ell} \mathrm{Gal}(E/F)$ above can sometimes be sharpened by a class field theoretic analysis similar to that of Lemma 3.2. Define

$$\Gamma_S = \prod_{v \in S_F} U_v/U_v^\ell,$$

where $S_F$ is the set of primes over $S$ in $F$. Recall that $\lambda$ denotes the prime over $\ell$ in $F$ and let $\mathcal{U}$ be the image in $\Gamma_S$ of those global units $\epsilon$ in $F$ such that $\epsilon \equiv 1 \bmod \lambda^2$. Then $\mathrm{Gal}(E/F)$ is a quotient of $\Gamma_S/\mathcal{U}$.

**Proposition 3.5.** *Let $A/\mathbb{Q}$ be a semistable abelian variety with good reduction at a prime $\ell \leq 19$ and set $H = \mathrm{Gal}(\mathbb{Q}(A[\ell])/\mathbb{Q}(\boldsymbol{\mu}_\ell))$. Then $\ell$ annihilates the maximal abelian quotient $H/H'$. If $H$ is nilpotent, then $H$ is an $\ell$-group.*

*Proof.* If $\ell$ is odd, or if $\ell = 2$ and $i \in \mathbb{Q}(A[\ell])$, let $L = \mathbb{Q}(A[\ell])$. Otherwise, let $L = \mathbb{Q}(i, A[2])$, in which case, $\mathrm{Gal}(L/\mathbb{Q}(i)) \simeq H$. Then $L$ is an $(\ell, S)$-controlled extension and our result follows from Lemma 3.4. Note that if $H$ is nilpotent, it is the product of all its Sylow subgroups, so it must be an $\ell$-group. $\square$

**Lemma 3.6.** *Suppose $L$ is an $(\ell, p)$-controlled extension containing $F = \mathbb{Q}(\boldsymbol{\mu}_{2\ell})$ and $\ell$ is a regular prime. Put $G = \mathrm{Gal}(L/\mathbb{Q})$ and assume $H = \mathrm{Gal}(L/F)$ is an $\ell$-group. Let $\tilde{\sigma}$ be a generator for the inertia group $\mathcal{I}_\mathfrak{P}$ of a prime $\mathfrak{P}$ over $p$ in $L$ and suppose $\tilde{\tau} \in G$ restricts to a generator of $\mathrm{Gal}(F/\mathbb{Q})$. Then the conjugates of $\tilde{\sigma}$ by powers of $\tilde{\tau}$ generate $H$, while $\tilde{\sigma}$ and $\tilde{\tau}$ generate $G$. If $p$ is unramified in $L$, then $L = F$.*

*Proof.* Note that $\mathcal{I}_\mathfrak{P}$ is cyclic of order dividing $\ell$ by (L3). Let $N$ be the subgroup of $H$ generated by $H'$, $H^\ell$ and the conjugates of $\tilde{\sigma}$ by powers of $\tilde{\tau}$. Then $N$ is normal in $G$ and corresponds to an $(\ell, p)$-controlled extension $M/F$, unramified outside $\ell$, such that $\mathrm{Gal}(M/F)$ is an elementary abelian $\ell$-group. So $N = H$ by the Remark after Lemma 3.2. Our claim now follows from Burnside's lemma. $\square$

Next, we give a description of the maximal $(2, p)$-controlled 2-extension of $\mathbb{Q}$. Since we do not need this result, we leave the proof to the reader as an exercise in class field theory. A general study of extensions of number fields with wild ramification of bounded depth can be found in recent work of Hajir and Maire [HM].

**Proposition 3.7.** *Consider $K = \mathbb{Q}(\sqrt{-p})$ and let $n$ be the 2-part of the class number of $K$. If $p \equiv 3 \bmod 4$, define $M = K(i)$. If $p \equiv 1 \bmod 4$, let $M$ be the 2-part of the ray class field of conductor 2 over $K$. Then $\mathrm{Gal}(M/K)$ is cyclic of order $2n$ and $\mathrm{Gal}(M/\mathbb{Q})$ is dihedral. Moreover, $M$ is $(2, p)$-controlled and contains every $(2, p)$-controlled 2-extension of $\mathbb{Q}$.*

Finally, we summarize for later use some elementary facts about pure submodules unrelated to the rest of this section.

**Lemma 3.8.** *Suppose $X$ and $Y$ are pure submodules of a free $\mathbb{Z}_\ell$-module $T$ of finite rank. Write $\bar{X} = (X + \ell T)/\ell T$ for the projection of $X$ to $T/\ell T$. Then:*
  (i) $\overline{X \cap Y} \subseteq \bar{X} \cap \bar{Y}$ *and* $\overline{X + Y} = \bar{X} + \bar{Y}$;
  (ii) $X \cap Y$ *is pure;*
  (iii) *if* $\bar{X} = 0$, *then* $X = 0$;
  (iv) *if* $\bar{X} \cap \bar{Y} = 0$, *then* $X + Y$ *is pure and a direct sum.*
*Let $T'$ be a free $\mathbb{Z}_\ell$-module of finite rank and suppose there is a perfect pairing $e: T \times T' \to \mathbb{T}_\ell(\boldsymbol{\mu})$. Then $X^\perp + Y^\perp \subseteq (X \cap Y)^\perp$, with equality if $X^\perp + Y^\perp$ is pure.*

*Proof.* Property (i) is clear and does not use purity. The natural injection

$$T/(X \cap Y) \hookrightarrow T/X \oplus T/Y$$

implies (ii). By Nakayama's Lemma and the isomorphism $\bar{X} \simeq X/(X \cap \ell T) = X/\ell X$ we have (iii). For (iv), suppose $\ell z = x + y$ for some $z \in T$, $x \in X$ and $y \in Y$. Then the coset $\bar{x} = -\bar{y}$ in $T/\ell T$ is an element of $\bar{X} \cap \bar{Y} = 0$. Hence $\bar{x} = \bar{y} = 0$. By purity of $X$ and $Y$, we may write $x = \ell x_1$ and $y = \ell y_1$ for some $x_1 \in X$ and $y_1 \in Y$. Since $T$ is torsion-free, we then have $z = x_1 + y_1 \in X + Y$. Hence $T/(X + Y)$ is torsion-free and (iv) is verified. Given the perfect pairing $e$, we clearly have $X^\perp + Y^\perp \subseteq (X \cap Y)^\perp$. Equality holds if $X^\perp + Y^\perp$ is pure because both sides have the same rank. □

## 4. Decomposition of the Tate module

Suppose $A/\mathbb{Q}$ is an abelian variety with good reduction at $\ell$ and semistable bad reduction at $p$. Grothendieck described certain submodules of the Tate module $\mathbb{T}_\ell(A)$ that are Galois modules for a decomposition group over $p$ in $G_\mathbb{Q}$. A convenient summary of the information we need also appears in [Ed, §2]. In this section, we use our assumptions that the $\ell$-division field of $A$ is small and that $A$ has good reduction outside $p$ to create $G_\mathbb{Q}$-submodules of $\mathbb{T}_\ell(A)$ from these Grothendieck modules. First we establish the relevant notation.

The connected component of the identity $\mathcal{A}_p^0$ of the special fiber $\mathcal{A}_p$ of the Néron model of $A$ at $p$ admits a decomposition $0 \to \mathcal{T} \to \mathcal{A}_p^0 \to \mathcal{B} \to 0$ in which $\mathcal{T}$ is a torus and $\mathcal{B}$ an abelian variety defined over $\mathbb{F}_p$. Setting $\dim \mathcal{T} = t$ and $\dim \mathcal{B} = a$, we have $t + a = \dim A$. Write $\Phi_A = \mathcal{A}_p/\mathcal{A}_p^0$ for the group of connected components. Denote the perfect pairing on the Tate modules of $A$ and its dual abelian variety $\hat{A}$ by

(4.1) $$e_\infty : \mathbb{T}_\ell(A) \times \mathbb{T}_\ell(\hat{A}) \to \mathbb{T}_\ell(\boldsymbol{\mu}) = \varprojlim \boldsymbol{\mu}_{\ell^n}.$$

Let $L_\infty = \mathbb{Q}(A[\ell^\infty])$ be the $\ell$-division tower of $A$ and set $G_\infty = \mathrm{Gal}(L_\infty/\mathbb{Q})$. Since we may be moving among abelian varieties $\mathbb{Q}$-isogenous to $A$, it is important to note that $L_\infty$ only depends on the $\mathbb{Q}$-isogeny class of $A$. Clearly $G_\mathbb{Q}$ acts on $\mathbb{T}_\ell(A)$ through $G_\infty$. Fix an embedding $\iota : L_\infty \to \bar{\mathbb{Q}}_p$ and let $v$ be the corresponding valuation. Write $\mathcal{D}_v \supseteq \mathcal{I}_v$ for its decomposition and inertia groups in $G_\infty$. For $g \in G_\infty$, let $gv$ denote the valuation corresponding to $\iota \circ g^{-1}$.

As in [Gro, §2.5], we define $\mathcal{M}_1 = \mathcal{M}_1(A, v)$ to be the submodule of $\mathbb{T}_\ell(A)$ fixed by $\mathcal{I}_v$ and $\mathcal{M}_2 = \mathcal{M}_2(A, v)$ to be the submodule of $\mathbb{T}_\ell(A)$ orthogonal to $\mathcal{M}_1(\hat{A}, v)$ under the $e_\infty$-pairing. Clearly $g \in G_\infty$ acts by $g(\mathcal{M}_j(A, v)) = \mathcal{M}_j(A, gv)$ for $j = 1, 2$. Thus $\mathcal{M}_1$ and $\mathcal{M}_2$ are $\mathcal{D}_v$-modules.

The Igusa-Grothendieck theorem [Gro, Thm. 2.5] asserts that if $A$ is semistable at $p$, then $\mathcal{M}_2 \subseteq \mathcal{M}_1$. Further, the successive quotients in the decomposition

(4.2) $$\mathbb{T}_\ell(A) \supset \mathcal{M}_1 \supseteq \mathcal{M}_2 \supset 0,$$

are torsion-free $\mathbb{Z}_\ell$-modules. Both $\mathcal{M}_1$ and $\mathcal{M}_2$ are modules for $\mathcal{D}_v/\mathcal{I}_v \simeq G_{\mathbb{F}_p}$. We may identify $\mathcal{M}_2 \simeq \mathbb{T}_\ell(\mathcal{T})$ and $\mathcal{M}_1/\mathcal{M}_2 \simeq \mathbb{T}_\ell(\mathcal{B})$, so the $\mathbb{Z}_\ell$-ranks of $\mathcal{M}_2$ and $\mathcal{M}_1/\mathcal{M}_2$ are $t$ and $2a$ respectively. Using the $e_\infty$-pairing, one can see that $(g - 1)(\mathbb{T}_\ell(A)) \subseteq \mathcal{M}_2$ for all $g \in \mathcal{I}_v$.

Let $\mathbb{T}_\ell(\varphi)$ denote the map of Tate modules induced by a $\mathbb{Q}$-isogeny $\varphi : A \to A'$ and set $\mathcal{M}'_j = \mathcal{M}_j(A', v)$ for $j = 1, 2$. Then $\mathbb{T}_\ell(\varphi)(\mathcal{M}_j) \subseteq \mathcal{M}'_j$. We can also show that $\mathbb{T}_\ell(\varphi)^{-1}(\mathcal{M}'_j) \subseteq \mathcal{M}_j$ by using the existence of a quasi-inverse isogeny $\varphi' : A' \to A$ such that $\varphi' \circ \varphi$ (resp. $\varphi \circ \varphi'$) is multiplication by $m$ on $A$ (resp. $A'$), where $m$ is the exponent of the kernel of $\varphi$. Indeed, if $\mathbb{T}_\ell(\varphi)(x) \in \mathcal{M}'_j$, then $mx \in \mathbb{T}_\ell(\varphi')(\mathcal{M}'_j) \subseteq \mathcal{M}_j$. But $\mathbb{T}_\ell(A)/\mathcal{M}_j$ is torsion-free, so $x \in \mathcal{M}_j$.

The following result, mildly strengthening [BK2, Lem. 3], describes a more subtle effect of isogeny. Write $\bar{\mathcal{M}}_1^{(n)}$ and $\bar{\mathcal{M}}_2^{(n)}$ for the projections of $\mathcal{M}_1$ and $\mathcal{M}_2$ to the $n^{\mathrm{th}}$ layer $\mathbb{T}_\ell(A)/\ell^n \mathbb{T}_\ell(A) \simeq A[\ell^n]$. Write $|X|_\ell$ for the order of the $\ell$-primary subgroup of a finite abelian group $X$. See (1.3) for our definition of $\ell$-maximal abelian varieties.

**Lemma 4.3.** *Suppose $A/\mathbb{Q}$ has good reduction at $\ell$ and semistable bad reduction at $p$. Let $\varphi : A \to A'$ be a $\mathbb{Q}$-isogeny whose kernel $\kappa$ has exponent $\ell^n$. Then*

$$|\Phi_{A'}(\bar{\mathbb{F}}_p)|_\ell \cdot |\kappa/(\kappa \cap \bar{\mathcal{M}}_1^{(n)})| = |\Phi_A(\bar{\mathbb{F}}_p)|_\ell \cdot |\kappa \cap \bar{\mathcal{M}}_2^{(n)}|.$$

*In particular, if $A$ is $\ell$-maximal and $\kappa$ is a Galois submodule of $A[\ell^n]$ contained in $\bar{\mathcal{M}}_1^{(n)}$, then $\kappa \cap \bar{\mathcal{M}}_2^{(n)} = 0$.*

*Proof.* The main argument is sketched in [BK2, Lem. 3] for $n = 1$ and carries over with obvious modifications. As a simplification to the formula given there, we may use the equality $|\Phi_A(\bar{\mathbb{F}}_p)| = |\Phi_{\hat{A}}(\bar{\mathbb{F}}_p)|$, arising from the existence of a pairing

$$\Phi_A(\bar{\mathbb{F}}_p) \times \Phi_{\hat{A}}(\bar{\mathbb{F}}_p) \to \mathbb{Q}/\mathbb{Z}, \tag{4.4}$$

defined by Grothendieck and shown to be non-degenerate for semistable abelian varieties by various authors. For this pairing and its history, see [Mc]. □

**Corollary 4.5.** *If $A$ has semistable bad reduction at an odd prime $p$ and is 2-maximal for $p$, then $\mathbb{Q}(A[2]) \neq \mathbb{Q}$. If, in addition, $\mathbb{Q}(A[2])$ is unramified outside 2 and $\infty$, then $\mathbb{Q}(A[2]) = \mathbb{Q}(i)$.*

*Proof.* If $\mathbb{Q}(A[2]) = \mathbb{Q}$, then $\bar{\mathcal{M}}_2^{(1)}$ is a non-trivial Galois module, in contradiction to the last assertion of the previous lemma. Under the additional assumption that $\mathbb{Q}(A[2])$ is unramified outside 2 and $\infty$, Lemma 3.6 implies that $L = \mathbb{Q}(A[2], i)$ is equal to $\mathbb{Q}(i)$. It follows that $\mathbb{Q}(A[2]) = \mathbb{Q}(i)$. □

Recall that there is a dictionary (see [Sch, §2.5]) between $G_\infty$-submodules $X$ of $T = \mathbb{T}_\ell(A)$ and abelian varieties $\mathbb{Q}$-isogenous to $A$. If $\ell^n T \subseteq X$ then $\kappa = X/\ell^n T \subseteq A[\ell^n]$ is the kernel of an isogeny $\varphi : A \to A' = A/\kappa$. Write $\varphi'$ for the quasi-inverse isogeny, such that $\varphi' \circ \varphi$ is multiplication by $\ell^n$ on $A$ and $\varphi \circ \varphi'$ is multiplication by $\ell^n$ on $A'$. Then $\mathbb{T}_\ell(\varphi')$ provides a Galois isomorphism of $\mathbb{T}_\ell(A')$ onto $X$. In particular, $X/\ell X \simeq \mathbb{T}_\ell(A')/\ell \mathbb{T}_\ell(A') \simeq A'[\ell]$.

**Lemma 4.6.** *Let $A/\mathbb{Q}$ be an abelian variety and let $W$ be a $\mathbb{Z}_\ell$-submodule of $T = \mathbb{T}_\ell(A)$ stabilized by a subgroup $\mathcal{H}$ of $G_\infty$. Define $W_n = W + \ell^n T$ and let $\mathcal{S}$ be the set of integers $n \geq 0$ such that $W_n$ is stabilized by $G_\infty$. Let $\mathcal{C}_0$ be a collection of abelian varieties containing at least those that correspond to the $G_\infty$-submodules $W_n$ for $n \in \mathcal{S}$. Define*

$$\mathcal{N} = \{g \in G_\infty \,|\, g \text{ acts by homothety on } A'[\ell] \text{ for all } A' \in \mathcal{C}_0\}.$$

*If $\mathcal{H}$ and $\mathcal{N}$ generate $G_\infty$, then $G_\infty$ stabilizes $W$.*

*Proof.* Since $W_0 = T$, we have $0 \in \mathcal{S}$. Suppose $n \in \mathcal{S}$. Then $W_n$ is stabilized by $G_\infty$, so $W_n$ corresponds to $A' = A/\kappa_n$, where $\kappa_n = W_n/\ell^n T$. By assumption, $A'$ is in $\mathcal{C}_0$, so the elements of $\mathcal{N}$ act as homotheties on $A'[\ell] \simeq W_n/\ell W_n$. For $w \in W$, consider a coset $w + \ell W_n$ in $W_n/\ell W_n$. If $g \in \mathcal{N}$ acts on this coset by multiplication by $\alpha$, then $g(w) + \ell W_n = \alpha w + \ell W_n$ and therefore

$$g(w) \in \alpha w + \ell W_n = \alpha w + \ell(W + \ell^n T) \subseteq W + \ell^{n+1} T = W_{n+1}.$$

Hence $\mathcal{N}$ stabilizes $W_{n+1}$. Given that $G_\infty$ is generated by $\mathcal{H}$ and $\mathcal{N}$ and that $\mathcal{H}$ already stabilizes $W$, we find that $G_\infty$ stabilizes $W_{n+1}$. Therefore $n + 1$ is in $\mathcal{S}$ and it follows that $G_\infty$ stabilizes $W_n$ for all $n \geq 0$. Hence $G_\infty$ stabilizes $W$. □

We now explain the special assumptions on which our study of abelian varieties with small $\ell$-division fields rests. Suppose $A/\mathbb{Q}$ is semistable, with bad reduction only at $p$ and that $H = \text{Gal}(\mathbb{Q}(A[\ell])/\mathbb{Q}(\boldsymbol{\mu}_\ell))$ is an $\ell$-group. Let $\mathcal{C}$ denote the $\mathbb{Q}$-isogeny class of $A$. Thanks to the theorem of Faltings [Fa, Satz 6] and its extension to unpolarized abelian varieties by Zarhin [Za, Thm. 1], the set of isomorphism classes in $\mathcal{C}$ is finite. Put $L$ for the compositum of the $\ell$-division fields of the varieties in $\mathcal{C}$ and, if necessary, adjoin $i = \sqrt{-1}$ when $\ell = 2$. Then $L$ contains $F = \mathbb{Q}(\boldsymbol{\mu}_{2\ell})$ and $L$ is an $(\ell, p)$-controlled extension, as defined in section 3.

Once and for all, we choose a topological generator $\sigma$ for the cyclic pro-$\ell$ group $\mathcal{I}_v \subset G_\infty$. There exists an element $\tau \in G_\infty$ whose restriction to $F$ generates $\text{Gal}(F/\mathbb{Q})$. When $\ell$ is odd, we choose $\tau$ to have order $\ell - 1$, using the fact that $H$ is an $\ell$-group. When $\ell = 2$, we choose $\tau$ to be the generator of some inertia group at an archimedean place of $L_\infty$, i.e. a complex conjugation. In view of Lemma 3.6, the restrictions of $\sigma$ and $\tau$ to $L$ generate $\text{Gal}(L/\mathbb{Q})$ when $\ell$ is a regular prime. If $\ell = 2$ or $3$ then $\tau^2 = 1$. However, our methods only require that the action of $\tau$ on $\mathbb{T}_\ell(A)$ satisfy a polynomial of degree at most 2, so they also apply, for example, if $A$ is of $\text{GL}_2$-type.

**Basic Assumptions.**

(C1) $A/\mathbb{Q}$ is semistable, with bad reduction only at $p$;
(C2) $\ell$ is a regular prime and $H = \text{Gal}(\mathbb{Q}(A[\ell])/\mathbb{Q}(\boldsymbol{\mu}_\ell))$ is an $\ell$-group;
(C3) the action of $\tau$ on $\mathbb{T}_\ell(A)$ satisfies a polynomial of degree at most 2.

These conditions depend only on the $\mathbb{Q}$-isogeny class $\mathcal{C}$ of $A$. Indeed, $\mathbb{Q}(B[\ell^\infty]) = L_\infty$ is the same for all $B$ in $\mathcal{C}$. Furthermore, $H$ is an $\ell$-group if and only if $\text{Gal}(\mathbb{Q}(B[\ell])/\mathbb{Q}(\boldsymbol{\mu}_\ell))$ is an $\ell$-group, since $\text{Gal}(L_\infty/\mathbb{Q}(B[\ell]))$ is pro-$\ell$. Finally, the minimal polynomial of $\tau$ can be read from $\mathbb{T}_\ell(A) \otimes \mathbb{Q}_\ell$.

We will show, under (C1) and (C2), that the following $\mathbb{Z}_\ell$-submodules of $\mathbb{T}_\ell(A)$ are in fact $G_\infty$-modules:

$$(4.7) \qquad W_A = \bigcap_{j \geq 0} \tau^j(\mathcal{M}_1(A, v)) \quad \text{and} \quad Y_A = \sum_{j \geq 0} \tau^j(\mathcal{M}_2(A, v)).$$

**Proposition 4.8.** *Suppose $A/\mathbb{Q}$ satisfies (C1) and (C2). Then $W_A$ and $Y_A$ are stabilized by $G_\infty$. Furthermore $W_A$ is a pure submodule of $\mathbb{T}_\ell(A)$.*

*Proof.* Let $\mathcal{H}$ be the closed subgroup of $G_\infty$ generated by $\sigma$ and $\tau$. Clearly, $\tau$ acts on $W_A$ and $Y_A$, in view of the definitions (4.7). Since $\sigma$ is the identity map on $\mathcal{M}_1(A, v)$, it acts as the identity on $W_A$. As noted above, $(\sigma - 1)(\mathbb{T}_\ell(A)) \subseteq \mathcal{M}_2(A, v)$, so $\sigma$ acts on $Y_A$. Hence $W_A$ and $Y_A$ are modules for $\mathcal{H}$.

Let $L$ be the compositum of the $\ell$-division fields of all the abelian varieties $\mathbb{Q}$-isogenous to $A$ and, if necessary when $\ell = 2$, also adjoin $i$. We wish to apply Lemma 4.6, with $T = \mathbb{T}_\ell(A)$, $\mathcal{G} = G_\infty$ and $W = W_A$ or $Y_A$. Then $\mathcal{N}$ certainly contains $\text{Gal}(L_\infty/L)$. According to Lemma 3.6, the restrictions of $\sigma$ and $\tau$ to $L$ generate $\text{Gal}(L/\mathbb{Q})$. Therefore $\mathcal{H}$ and $\mathcal{N}$ generate $G_\infty$. It follows that $W_A$ and $Y_A$ are $G_\infty$-modules.

The purity of $\mathcal{M}_1(A, v)$ implies the purity of $W_A$ as a $\mathbb{Z}_\ell$-submodule of $\mathbb{T}_\ell(A)$ by Lemma 3.8. □

Recall from (1.3) the definition of an $\ell$-maximal abelian variety. If $X$ is a $\mathbb{Z}_\ell$-submodule of $\mathbb{T}_\ell(A)$, we write $\bar{X} = (X + \ell\mathbb{T}_\ell(A))/\ell\mathbb{T}_\ell(A)$ for the projection of $X$ to $\mathbb{T}_\ell(A))/\ell\mathbb{T}_\ell(A) \simeq A[\ell]$.

**Lemma 4.9.** *Let $A$ be an abelian variety satisfying (C1) and (C2). Put $\bar{\mathcal{M}}_1 = \bar{\mathcal{M}}_1(A, v)$ and define*
$$\kappa = \bigcap_{j \geq 0} \tau^j(\bar{\mathcal{M}}_1).$$
*If $A$ is $\ell$-maximal, then $\kappa \cap \bar{\mathcal{M}}_2 = 0$.*

*Proof.* By Lemma 3.6, the restrictions of $\sigma$ and $\tau$ generate $\mathrm{Gal}(\mathbb{Q}(A[\ell])/\mathbb{Q})$. But $\sigma$ acts trivially on $\mathcal{M}_1$ and therefore on $\kappa$, while $\tau$ acts on $\kappa$ from the definition. Hence $\kappa$ is a Galois invariant subspace of $A[\ell]$ contained in $\bar{\mathcal{M}}_1$ and we conclude by Lemma 4.3. □

**Corollary 4.10.** *If $A$ satisfies (C1) and (C2), then there are at least 2 primes over $p$ in $F = \mathbb{Q}(\boldsymbol{\mu}_{2\ell})$. In addition, if $\ell = 2$, then $p \equiv 1 \bmod 8$.*

*Proof.* We may assume that $A$ is $\ell$-maximal. Suppose $\ell$ is odd and there is only one prime over $p$ in $F$. Then the decomposition group $\mathcal{D}_v$ projects onto the group $\mathrm{Gal}(F/\mathbb{Q})$ of order $\ell - 1$, with pro-$\ell$ kernel. Hence we could have chosen $\tau$ in $\mathcal{D}_v$ and so $\kappa$ in the lemma is simply $\bar{\mathcal{M}}_1$. Therefore $\bar{\mathcal{M}}_2 = \kappa \cap \bar{\mathcal{M}}_2 = 0$. This implies that the toroidal dimension of the bad fiber of $A$ is 0, contradicting bad reduction at $p$.

For $\ell = 2$, we showed in [BK2, Prop. 5] that $p \equiv 1 \bmod 4$, whence there are 2 primes over $p$ in $F = \mathbb{Q}(i)$. Suppose $p \equiv 5 \bmod 8$ and let $\phi$ be a Frobenius element of $\mathcal{D}_v$, in the sense that $\phi$ induces the $p^{\mathrm{th}}$ power map on the residue field. Then $\mathcal{D}_v = \langle \phi, \sigma \rangle$ and $\phi\sigma\phi^{-1} = \sigma^p$. Since $\phi(\sqrt{2}) = -\sqrt{2}$, the restrictions of $\tau, \phi$ and $\sigma$ to $E = \mathbb{Q}(i, \sqrt{2}, \sqrt{p})$ clearly generate $\mathrm{Gal}(E/\mathbb{Q})$. But $E$ is the maximal elementary 2-extension of $\mathbb{Q}$ unramified outside $\{2, p, \infty\}$, so $G_\infty = \langle \tau, \phi, \sigma \rangle$ by Burnside's Lemma. Hence $\mathcal{D}_v$ has index 2, so is normal in $G_\infty$.

Now we may write $\tau\sigma\tau^{-1} = \phi^x \sigma^y$ for some $x, y \in \mathbb{Z}_2$. Checking the action of both sides of this equation on a generator for $\mathbb{T}_2(\boldsymbol{\mu})$, we find that $x = 0$, so $\tau$ normalizes $\mathcal{I}_v$. Hence $\mathcal{M}_1$ is stabilized by $\tau$ and we get a contradiction as in the previous case. Therefore $p \equiv 1 \bmod 8$. □

From hereon in, we impose the additional assumption (C3). Since the dual abelian variety $\hat{A}$ is isogenous to $A$, it also satisfies (C3). Write $\bar{\mathcal{M}}'_j$ for the projection of $\mathcal{M}'_j = \mathcal{M}_j(\hat{A}, v)$ to $\hat{A}[\ell]$. In the perfect pairing

(4.11) $$A[\ell] \times \hat{A}[\ell] \to \boldsymbol{\mu}_\ell$$

induced by the $e_\infty$-pairing (4.1), we have $(\bar{\mathcal{M}}_1)^\perp = \bar{\mathcal{M}}'_2$. Recall that $t$ and $a$ respectively denote the toroidal and abelian dimension of the bad fiber of $A$ at $p$.

**Proposition 4.12.** *Suppose $A$ satisfies (C1), (C2), (C3) and is $\ell$-maximal. Then $W_A$ and $Y_A$ are pure $\mathbb{Z}_\ell$-submodule of $\mathbb{T}_\ell(A)$ of rank $2a$ and $2t$ respectively, stabilized by $G_\infty$. Under the $e_\infty$-pairing, $W_A^\perp = Y_{\hat{A}}$ and $Y_A^\perp = W_{\hat{A}}$. We have a direct sum decomposition $\mathbb{T}_\ell(A) = W_A \oplus Y_A$. Furthermore, the natural map $W_A \to \mathcal{M}_1(A, v)/\mathcal{M}_2(A, v)$ is an isomorphism of $\mathbb{Z}_\ell$-module.*

*Proof.* Hypothesis (C3) implies that $\tau^2(\mathcal{M}_2) \subseteq \mathcal{M}_2 + \tau(\mathcal{M}_2)$, so
$$Y_A = \sum_{j \geq 0} \tau^j(\mathcal{M}_2) = \mathcal{M}_2 + \tau(\mathcal{M}_2).$$

Furthermore, $\tau^2(\bar{\mathcal{M}}'_2) \subseteq \bar{\mathcal{M}}'_2 + \tau(\bar{\mathcal{M}}'_2)$ and so $\tau^2(\bar{\mathcal{M}}_1) \supseteq \bar{\mathcal{M}}_1 \cap \tau(\bar{\mathcal{M}}_1)$ by the $\mathbb{F}_\ell$-vector space duality (4.11). Hence the Galois invariant subspace $\kappa$ of Lemma 4.9 reduces to $\kappa = \bar{\mathcal{M}}_1 \cap \tau(\bar{\mathcal{M}}_1)$ and $\kappa \cap \bar{\mathcal{M}}_2 = 0$ by that Lemma. Since $\mathcal{M}_2 \subseteq \mathcal{M}_1$, we have

$$\tau(\bar{\mathcal{M}}_2) \cap \bar{\mathcal{M}}_2 \subseteq \tau(\bar{\mathcal{M}}_1) \cap \bar{\mathcal{M}}_2 = \kappa \cap \bar{\mathcal{M}}_2 = 0. \tag{4.13}$$

Therefore, $Y_A = \mathcal{M}_2 + \tau(\mathcal{M}_2)$ is pure and a direct sum, of rank $2t$, by Lemma 3.8(iv). Similarly, $Y_{\hat{A}}$ is pure in $\mathbb{T}_\ell(\hat{A})$, so we have

$$Y_{\hat{A}} = \sum_{j \geq 0} \tau^j(\mathcal{M}'_2) = \sum_{j \geq 0} \tau^j(\mathcal{M}_1^\perp) = (\bigcap_{j \geq 0} \tau^j(\mathcal{M}_1))^\perp = W_A^\perp.$$

Hence $W_A$ is pure of rank $2a$. In view of Lemma 3.8(i) and (4.13), we have

$$\bar{W}_A \cap \bar{\mathcal{M}}_2 \subseteq (\bar{\mathcal{M}}_1 \cap \tau(\bar{\mathcal{M}}_1)) \cap \bar{\mathcal{M}}_2 = \tau(\bar{\mathcal{M}}_1) \cap \bar{\mathcal{M}}_2 = 0.$$

Hence $W_A + \mathcal{M}_2$ is pure of rank $t + 2a$ and a direct sum. By the obvious inclusion and equality of ranks, we have $W_A + \mathcal{M}_2 = \mathcal{M}_1$. It follows that

$$W_A \simeq (W_A + \mathcal{M}_2)/\mathcal{M}_2 = \mathcal{M}_1/\mathcal{M}_2.$$

By (4.13) and Lemma 3.8(iv), we also find that $\tau(\mathcal{M}_1) + \mathcal{M}_2$ is a pure submodule of $\mathbb{T}_\ell(A)$ and a direct sum. Equality of ranks implies that $\mathbb{T}_\ell(A) = \tau(\mathcal{M}_1) \oplus \mathcal{M}_2$ and so $\mathbb{T}_\ell(A) = \tau(W_A) \oplus \tau(\mathcal{M}_2) \oplus \mathcal{M}_2 = W_A \oplus Y_A$, as claimed. Finally, we have already observed in Proposition 4.8 that $W_A$ and $Y_A$ are $G_\infty$-modules. $\square$

## 5. The $\ell$-adic formal group and the proof of theorem 1.1

Suppose that $A$ is an abelian variety of dimension $d$ defined over $\mathbb{Q}_\ell$, with good reduction modulo $\ell$. The kernel of reduction has the structure of a formal group of height $h$ and dimension $d$, with $d \leq h \leq 2d$. We put $\mathcal{F}_A$ for this formal group, suppressing the dependence on $\ell$. The $G_{\mathbb{Q}_\ell}$-module $C = \mathcal{F}_A[\ell^\infty]$ is isomorphic to $(\mathbb{Q}_\ell/\mathbb{Z}_\ell)^h$ as an abelian group. Furthermore, $C$ is a connected $\ell$-divisible group over $\mathbb{Z}_\ell$ while $A[\ell^\infty]/C$ is étale.

Set $\mathbb{V}_\ell(A) = \mathbb{T}_\ell(A) \otimes_{\mathbb{Z}_\ell} \mathbb{Q}_\ell$ and $\mathbb{V}_\ell(\mathcal{F}_A) = \mathbb{T}_\ell(\mathcal{F}_A) \otimes_{\mathbb{Z}_\ell} \mathbb{Q}_\ell$ for the Tate vector spaces of $A$ and $\mathcal{F}_A$. Clearly the inertia group inside $G_{\mathbb{Q}_\ell}$ acts trivially on $\mathbb{V}_\ell(A)/\mathbb{V}_\ell(\mathcal{F}_A)$, as this quotient corresponds to the Tate vector space of the reduction of $A$ over $\mathbb{F}_\ell$. In the special case of *ordinary reduction* (*i.e.* $h = d$), we have

$$\mathbb{V}_\ell(\mathcal{F}_A)^\perp = \mathbb{V}_\ell(\mathcal{F}_{\hat{A}}) \tag{5.1}$$

under the Weil pairing $\mathbb{V}_\ell(A) \times \mathbb{V}_\ell(\hat{A}) \to \mathbb{Q}_\ell(1)$, where $\hat{A}$ denotes the dual abelian variety of $A$. (See for example [CG, p. 154].)

**Lemma 5.2.** *Let $A$ be an abelian variety of dimension $d$ defined over $\mathbb{Q}_\ell$ with good reduction modulo $\ell$. If the tame ramification degree of the extension $\mathbb{Q}_\ell(A[\ell])/\mathbb{Q}_\ell$ divides $\ell - 1$, then $A$ is ordinary.*

*Proof.* We extend the base to the ring of integers $R$ in the unramified closure $K = \mathbb{Q}_\ell^{\mathrm{nr}}$ to obtain an algebraically closed residue field. According to Raynaud (see [Gru, Thm. 4.4]), the finite group scheme $A[\ell]$ admits a composition series whose simple constituents are $\mathbb{F}$-vector space schemes for varying finite fields $\mathbb{F}$ of characteristic $\ell$. Such a constituent $M$ is a cyclic $\mathbb{F}$-module, determined as follows. By simplicity, the wild ramification subgroup of $G_K$ must act trivially on $M$, so the action of $G_K$ factors through a quotient $\mathrm{Gal}(E/K)$, such that $E \subseteq K(A[\ell])$ and $E/K$ is tamely ramified. Suppose $[E : K] = r$. Then $\mathbb{F} = \mathbb{F}_\ell(\boldsymbol{\mu}_r)$, with the action of a primitive $r^{\mathrm{th}}$ root of unity induced by the action of a generator for $\mathrm{Gal}(E/K)$.

In our case, $r$ divides $\ell - 1$, whence $\mathbb{F} = \mathbb{F}_\ell$ and $M$ therefore is a group scheme of order $\ell$. These are classified by Oort-Tate. Thus each simple constituent $M_j$ of $A[\ell]$ has the form $M_j \simeq G^\ell_{a_j, b_j}$, with $a_j b_j = \ell$ in the standard notation [Gru, Theorem 2.1]. Since the valuation $v$ of $R$ is unramified, we have $v(a_j) \in \{0, 1\}$. Then $G^\ell_{a_j, b_j}$ is étale (resp. connected) if $v(a_j) = 0$ (resp. 1). Let $n_0$ (resp. $n_1$) be the number of $a_j$ such that $v(a_j) = 0$ (resp. 1).

The exponent of the different of an affine group scheme $\mathcal{G}$ may be defined as follows [Gru, p. 62]. Suppose $\mathcal{G} = \mathrm{spec}(\mathfrak{A})$, where $\mathfrak{A}$ is an $R$-algebra. Let $\Omega^1_{\mathfrak{A}/R}$ be the module of Kähler differentials of $\mathfrak{A}$ and define $\delta(\mathcal{G})$ to be the length of the $R$-module $s^* \Omega^1_{\mathfrak{A}/R} = \Omega^1_{\mathfrak{A}/R} \otimes_\mathfrak{A} R$, where the tensor product is formed over the zero-section $s : \mathfrak{A} \to R$.

Since the $\ell$-divisible group arising from $\ell$-power torsion on $A$ has dimension $d$, we find that $\delta(A[\ell]) = d\, v(\ell) = d$ for example by [Gru, Prop. 3.4]. Furthermore, $\delta(G^\ell_{a_j, b_j}) = v(a_j)$ by [Gru, Prop. 2.3]. But $\delta$ behaves like an additive Euler characteristic on short exact sequences. Hence $d = \sum_j v(a_j) = n_1$. Note that the simple constituents that belong to the connected group scheme $\mathcal{F}_A[\ell]$ are precisely those for which $v(a_j) = 1$. By definition of the height $h$ of $\mathcal{F}_A$ we have

$$\ell^h = |\mathcal{F}_A[\ell]| = \prod_{v(a_j) = 1} |G^\ell_{a_j, b_j}| = \ell^{n_1}.$$

Hence $h = n_1 = d$ as claimed. □

We now return to the global situation: $A$ is an abelian variety defined over $\mathbb{Q}$ and we impose the basic assumptions (C1) and (C2) of section 4. Preserving previous notation, we set $L_\infty = \mathbb{Q}(A[\ell^\infty])$ and $G_\infty = \mathrm{Gal}(L_\infty/\mathbb{Q})$. As usual, we have a topological generator $\sigma$ for the inertia group $\mathcal{I}_v$ inside $G_\infty$ at a place $v$ over $p$ and we have an element $\tau \in G_\infty$ of order $[F : \mathbb{Q}]$, where $F = \mathbb{Q}(\boldsymbol{\mu}_{2\ell})$, such that the restriction of $\tau$ generates $\mathrm{Gal}(F/\mathbb{Q})$. According to Lemma 3.6, the restrictions of $\sigma$ and $\tau$ to the $\ell$-division field of $A$ generate its Galois group. By fixing an embedding of $L_\infty$ to $\bar{\mathbb{Q}}_\ell$, we get a natural $\mathcal{D}_v$-isomorphism $\mathbb{T}_\ell(A) \to \mathbb{T}_\ell(A_{\mathbb{Q}_\ell})$, where $A_{\mathbb{Q}_\ell}$ is obtained by base extension to $\mathbb{Q}_\ell$. We are interested in elements of $\mathbb{T}_\ell(A)$ whose image lands in $\mathbb{T}_\ell(\mathcal{F}_A) \subseteq \mathbb{T}_\ell(A_{\mathbb{Q}_\ell})$.

**Lemma 5.3.** *Let $A/\mathbb{Q}$ be an abelian variety satisfying (C1) and (C2). Suppose $X$ is a $G_\infty$-submodule of $\mathbb{T}_\ell(A)$ with trivial action by the inertia group $\mathcal{I}_v$. Let $K_\infty$ be the fixed field of the kernel of the representation of $G_\infty$ afforded by $X$. Then $K_\infty$ is unramified outside $\ell$ and totally ramified over $\ell$. Furthermore $X \subseteq \mathbb{T}_\ell(\mathcal{F}_A)$.*

*Proof.* By assumption, $K_\infty$ is unramified at $v$. But then $K_\infty/\mathbb{Q}$ is unramified at all places over $p$ because it is Galois.

It is convenient to pass to the compositum $K_\infty F$. Let $\lambda$ be the place of $K_\infty F$ over $\ell$ determined by our choice of embedding $L_\infty \to \bar{\mathbb{Q}}_\ell$. Since the prime over $\ell$ in $F$ is totally ramified, we may verify that $K_\infty$ is totally ramified over $\ell$ by showing that the inertia group $\mathcal{I}_{F, \lambda}$ inside $H_\infty = \mathrm{Gal}(K_\infty F/F)$ is equal to $H_\infty$. Suppose on the contrary that $\mathcal{I}_{F, \lambda}$ is a proper subgroup of $H_\infty$. Condition (C2) implies that $H_\infty$ is pro-$\ell$, so there exists a subgroup of index $\ell$ in $H_\infty$, containing $\mathcal{I}_{F, \lambda}$. The corresponding fixed field is a cyclic extension of $F$ of degee $\ell$ unramified everywhere, whose existence contradicts our assumption that $\ell$ is a regular prime.

Put $\bar{X}^{(n)} = (X + \ell^n \mathbb{T}_\ell(A))/\ell^n \mathbb{T}_\ell(A) \subseteq A[\ell^n]$ for the $n^{\mathrm{th}}$ layer of $X$. Then we may represent elements of $X$ in the form $x = \varprojlim x_n$ with $x_n \in \bar{X}^{(n)}$ compatible under the projection maps. Since $\ell$ is totally ramified in $K_\infty$, the reduction $\bar{x}_n$ modulo $\lambda$ is an element of $\bar{A}(\mathbb{F}_\ell)$. Suppose the $\ell$-part of the order of $\bar{A}(\mathbb{F}_\ell)$ is $\ell^m$. For all $n \geq m+1$, it follows that $\bar{x}_{n-m} = \ell^m \bar{x}_n = 0$. Hence $x_{n'}$ is in the kernel of reduction for all layers $n' = n - m \geq 1$ and we have $X \subseteq \mathbb{T}_\ell(\mathcal{F}_A)$, as claimed. □

**Theorem 5.4.** *Let $\ell$ be a regular prime and let $A/\mathbb{Q}$ be a semistable abelian variety with bad reduction only at $p$, such that $H = \text{Gal}(\mathbb{Q}(A[\ell])/\mathbb{Q}(\boldsymbol{\mu}_\ell))$ is an $\ell$-group. Assume that the action of $\tau$ on $\mathbb{T}_\ell(A)$ satisfies a polynomial of degree at most 2. Then $A$ has ordinary reduction modulo $\ell$ and totally toroidal reduction modulo $p$.*

*Proof.* We may assume that among the members of its $\mathbb{Q}$-isogeny class $\mathcal{C}$, the variety $A$ is $\ell$-maximal; that is, $A \in \mathcal{V}_\ell(\mathcal{C})$ in the notation of (1.3). The dual abelian varieties $A$ and $\hat{A}$ are $\mathbb{Q}$-isogenous and $\hat{A}$ also is $\ell$-maximal because $|\Phi_{\hat{A}}(\bar{\mathbb{F}}_p)| = |\Phi_A(\bar{\mathbb{F}}_p)|$ by the perfect pairing (4.4). In view of Proposition 4.12, we have pure submodules $W_A$ of $\mathbb{T}_\ell(A)$ and $W_{\hat{A}}$ of $\mathbb{T}_\ell(\hat{A})$, of rank $2a$ and stabilized by $G_\infty$. Furthermore, the Weil pairing induces a perfect pairing

$$(5.5) \qquad W_A \times W_{\hat{A}} \to \mathbb{T}_\ell(\boldsymbol{\mu}).$$

Because $\mathcal{I}_v$ acts trivially on $W_A$, we find that $W_A \subseteq \mathbb{T}_\ell(\mathcal{F}_A)$ by Lemma 5.3. Similarly, $W_{\hat{A}} \subseteq \mathbb{T}_\ell(\mathcal{F}_{\hat{A}})$. Under the assumption that $H$ is an $\ell$-group, Lemma 5.2 implies that $A$ has ordinary reduction modulo $\ell$, so the Weil pairing on $W_A \times W_{\hat{A}}$ is trivial by (5.1). This contradicts the perfect pairing (5.5) unless $W_A = W_{\hat{A}} = 0$. Hence $a = 0$ and $A$ is totally toroidal at $p$. □

**Proof of Theorem 1.1.** Assume that $A/\mathbb{Q}$ is semistable, with bad reduction only at $p$, that $\ell = 2$ or 3 and that $H = \text{Gal}(\mathbb{Q}(A[\ell])/\mathbb{Q}(\boldsymbol{\mu}_\ell))$ is nilpotent. Then $H$ is in fact an $\ell$-group by Proposition 3.5. We have $p \equiv 1 \bmod 8$ when $\ell = 2$ and $p \equiv 1 \bmod 3$ when $\ell = 3$ by Corollary 4.10. Since $\tau^2 = 1$ when $\ell = 2$ or 3, we may conclude by Theorem 5.4. □

## 6. Small $\ell$-division fields

As building blocks for semistable abelian varieties with small $\ell$-division fields and bad reduction at only one prime $p$, we consider elliptic curves $E/\mathbb{Q}$ of conductor $p$ such that the mod-$\ell$ representation $\rho_{E,\ell} : G_\mathbb{Q} \to \text{Aut}(E[\ell]) \simeq \text{GL}_2(\mathbb{F}_\ell)$ is not surjective. Let $\Delta_E$ denote the minimal discriminant of $E$ and recall that $\text{ord}_p(\Delta_E) = -\text{ord}_p(j_E)$, where $j_E$ is the $j$-invariant. If $\ell \geq 7$ or $\text{ord}_p(j_E) \not\equiv 0 \bmod \ell$, then $E$ admits a $\mathbb{Q}$-isogeny of degree $\ell$ by [Se2, Prop. 21]. However, if $\ell \leq 7$ and $|\Delta_E|$ is an $\ell^{\text{th}}$ power, $E$ admits an isogeny of degree $\ell$, for example by [BK1, Prop. 9.2]. It then follows, as in [Se2, p. 307], that the isogeny class of $E$ contains a curve with a rational point of order $\ell$. The only such examples for odd $\ell$ occur in the well-known cases $(\ell, p) = (3, 19), (3, 37)$ and $(5, 11)$, referred to below as Miyawaki curves [Mi]. For $\ell = 2$, according to Neumann [Ne] or Setzer [Sz], we have $p = 17$ or $p = u^2 + 64$. In the latter case, there is a unique isogeny class consisting of the curves

$$(6.1) \qquad y^2 + xy = \begin{cases} x^3 + \dfrac{u-1}{4}x^2 - x, & \Delta = p, \\ x^3 + \dfrac{u-1}{4}x^2 + 4x + u, & \Delta = -p^2, \end{cases}$$

where the sign of $u$ is chosen to guarantee that $u \equiv 1 \bmod 4$.

Since the order of the connected component group of a semistable elliptic curve $C/\mathbb{Q}_p$ is $|\Phi_C(\overline{\mathbb{F}}_p)| = \operatorname{ord}_p(\Delta_C)$, one easily verifies that each Neumann-Setzer or Miyawaki isogeny class contains a unique curve $E_\ddagger$ with maximal connected component group. Its discriminant satisfies

$$(6.2) \qquad \operatorname{ord}_p(\Delta_\ddagger) = \begin{cases} 4 & \text{if } (\ell, p) = (2, 17), \\ \ell & \text{otherwise.} \end{cases}$$

From the Tate parametrization over $\mathbb{Q}_p$, we find that $\mathbb{Q}(E_\ddagger[\ell])$ is unramified over $p$ and so $\mathbb{Q}(E_\ddagger[\ell]) = F = \mathbb{Q}(\boldsymbol{\mu}_{2\ell})$ by Lemma 3.6 and Corollary 4.5. According to [MO], $E_\ddagger$ is the strong Weil curve in its isogeny class.

Write $E = E_\ddagger$, suppressing the dependence on $\ell$ and $p$, which should be clear in context. Recall that the $\ell^\infty$-division tower $L_\infty$ is the same for all curves in the isogeny class of $E$. Let $G_\infty = \operatorname{Gal}(L_\infty/\mathbb{Q})$ and observe that $\operatorname{Gal}(L_\infty/F)$ is pro-$\ell$. As usual, $\sigma \in G_\infty$ denotes a topological generator of the inertia group at a place $v$ over the bad prime $p$ and $\tau \in G_\infty$ denotes an element of order $[F : \mathbb{Q}]$ whose restriction to $F$ generates $\operatorname{Gal}(F/\mathbb{Q})$.

The Grothendieck module $\mathcal{M}_1(E, v) = \mathcal{M}_2(E, v)$ is a pure $\mathbb{Z}_\ell$-module of rank 1 in $\mathbb{T}_\ell(E)$ and we have $\mathbb{T}_\ell(E) = Y_E = \mathcal{M}_2(E, v) \oplus \tau(\mathcal{M}_2(E, v))$ by Proposition 4.12. Choose a generator $P$ for $\mathcal{M}_2(E, v)$. With respect to the generating set $\{P, \tau(P)\}$ for $\mathbb{T}_\ell(E)$, we obtain a matrix representation

$$\rho_E : G_\infty \to \operatorname{Aut}(\mathbb{T}_\ell(E)) \simeq \operatorname{GL}_2(\mathbb{Z}_\ell).$$

The Tate parametrization of $E$ over $\mathbb{Q}_p$ shows that

$$(6.3) \qquad \rho_E(\sigma) = \begin{pmatrix} 1 & s \\ 0 & 1 \end{pmatrix},$$

where $s$ equals $\operatorname{ord}_p(\Delta_\ddagger)$ up to multiplication by a unit in $\mathbb{Z}_\ell$, so $s \equiv 0 \bmod \ell$. Since the eigenvalues of $\tau$ are 1 and $\omega = \chi(\tau)$, where $\chi$ is the cyclotomic character, we have

$$(6.4) \qquad \rho_E(\tau) = \begin{pmatrix} 0 & -\omega \\ 1 & 1+\omega \end{pmatrix}.$$

The following lemma will be used to study abelian varieties $\mathbb{Q}$-isogenous to products of $E_\ddagger$ for fixed $\ell$ and $p$. It is sufficient for our purposes to state it for abelian varieties $A, B$ defined over $\mathbb{Q}$. The representation $\rho_B$ of $\operatorname{Gal}(\mathbb{Q}(B[\ell^\infty])/\mathbb{Q})$ afforded by $\mathbb{T}_\ell(B)$ naturally extends to the completed group ring

$$\Lambda_\ell(B) = \mathbb{Z}_\ell[[\operatorname{Gal}(\mathbb{Q}(B[\ell^\infty])/\mathbb{Q})]].$$

**Lemma 6.5.** *Let $\varphi : B \to A$ be a $\mathbb{Q}$-isogeny of abelian varieties. For given $g_1, g_2$ in $\operatorname{Gal}(\mathbb{Q}(B[\ell^\infty])/\mathbb{Q})$, assume there exists $\gamma \in \Lambda_\ell(B)$ such that $\rho_B(g_1) = \rho_B(g_2 + \ell \gamma)$. Then the restrictions of $g_1$ and $g_2$ to $A[\ell]$ are equal.*

*Proof.* Given $a \in A[\ell]$, we may find $b \in B$ of $\ell$-power order, say $\ell^n$, such that $\varphi(b) = a$. Since $g_1$ acts on $b$ through the representation $\rho_B$ modulo $\ell^n$, we have $g_1(b) = g_2(b) + \ell \gamma(b)$. Then

$$g_1(a) = g_1(\varphi(b)) = \varphi(g_1(b)) = \varphi(g_2(b) + \ell \gamma(b)) = g_2(a) + \gamma(\ell a) = g_2(a),$$

because $\varphi$ commutes with Galois and $\ell a = 0$. $\square$

**Proposition 6.6.** *Suppose the abelian variety $A$ is $\mathbb{Q}$-isogenous to a product of Neumann-Setzer or Miyawaki curves for fixed $\ell$ and $p$. Let $G = \operatorname{Gal}(\mathbb{Q}(A[\ell])/\mathbb{Q})$.*

(a) If $\ell = 2$, then $G$ is a quotient of $\mathbb{Z}/2 \times \mathbb{Z}/2$.
(b) If $\ell = 3$, then $G$ is a quotient of $\mathcal{S}_3$.
(c) If $\ell = 5$, then $H = \mathrm{Gal}(\mathbb{Q}(A[5])/\mathbb{Q}(\boldsymbol{\mu}_5))$ is abelian of exponent dividing 5 and rank at most 2. Furthermore, $\tilde{\tau}^{-2} h \tilde{\tau}^2 = h^{-1}$ for all $h \in H$.

*Proof.* Since $A$ is isogenous to $E^d$, where $d = \dim A$ and $E = E_\ddagger$, we have $\mathbb{Q}(A[\ell^\infty]) = \mathbb{Q}(E^d[\ell^\infty]) = L_\infty$. Furthermore, $\mathbb{T}_\ell(E^d) \simeq \oplus_1^d \mathbb{T}_\ell(E)$ affords the representation $\rho_{E^d} = \oplus_1^d \rho_E$.

Suppose $\ell = 2$ or $3$, so that $\rho_E(\tau) = \begin{pmatrix} 0 & 1 \\ 1 & 0 \end{pmatrix}$. Then $\rho_E(\sigma\tau - \tau\sigma^{-1}) = s\,\mathbf{1}$, where $\mathbf{1}$ is the identity in $M_2(\mathbb{Z}_\ell)$. Therefore $\rho_{E^d}(\sigma\tau - \tau\sigma^{-1}) = s\,\mathbf{1}$, where $\mathbf{1}$ now denotes the identity in $M_{2d}(\mathbb{Z}_\ell)$. By Lemma 6.5, the restrictions $\tilde{\sigma}$ and $\tilde{\tau}$ of $\sigma$ and $\tau$ to $G$ therefore satisfy $\tilde{\sigma}\tilde{\tau} = \tilde{\tau}\tilde{\sigma}^{-1}$. But $\tilde{\sigma}$ and $\tilde{\tau}$ generate $G$ by Lemma 3.6 and $\tilde{\sigma}^\ell = 1$ by property (L3) of $(\ell, S)$-controlled extensions. Hence $G$ is a quotient of $\mathbb{Z}/2 \times Z/2$ if $\ell = 2$ or a quotient of $\mathcal{S}_3$ if $\ell = 3$.

For $\ell = 5$, we find that $\rho_E(\sigma\tau^2 - \tau^2\sigma^{-1}) = (1+\omega)s\,\mathbf{1} \in M_2(\mathbb{Z}_5)$, so
$$\rho_{E^d}(\sigma\tau^2 - \tau^2\sigma^{-1}) = (1+\omega)s\,\mathbf{1} \in M_{2d}(\mathbb{Z}_5).$$
Now $\tilde{\sigma}\tilde{\tau}^2 = \tilde{\tau}^2\tilde{\sigma}^{-1}$ by Lemma 6.5. Acccording to Lemma 3.6, the subgroup $H$ of $G$ is generated by the conjugates of $\tilde{\sigma}$ under the action of $\tilde{\tau}$. Since $\tilde{\tau}^2$ normalizes $\tilde{\sigma}$, to prove that $H$ is abelian, it suffices to show that $\tilde{\sigma}$ commutes with $\tilde{\sigma}^\tau = \tilde{\tau}^{-1}\tilde{\sigma}\tilde{\tau}$. The reader may verify by judicious use of $\rho_E(\sigma - 1)$ and $\rho_E(\tau)$ that the image of $\rho_E$ on the group ring $\mathbb{Z}_5[[G_\infty]]$ contains all $2 \times 2$ matrices congruent to 0 modulo $s$. One can also check that
$$\rho_E(\sigma^\tau \sigma - \sigma\,\sigma^\tau) = -s^2\,\omega \begin{pmatrix} 1 & 2(1+\omega) \\ 0 & -1 \end{pmatrix}$$
Hence $\rho_E(\sigma^\tau \sigma - \sigma\,\sigma^\tau)$ belongs to $5\,\rho_E(\mathbb{Z}_5[[G_\infty]])$ and so $\rho_{E^d}(\sigma^\tau \sigma - \sigma\,\sigma^\tau)$ belongs to $5\,\rho_{E^d}(\mathbb{Z}_5[[G_\infty]])$. It follows from Lemma 6.5 that $\tilde{\sigma}^\tau$ and $\tilde{\sigma}$ do commute. We may conclude that $H$ is abelian of exponent dividing 5 and rank at most 2, generated by $\tilde{\sigma}$ and $\tilde{\sigma}^\tau$. Furthermore, the action of conjugation by $\tilde{\tau}^2$ on $H$ is given by inversion. $\square$

**Proposition 6.7.** *Let $\mathcal{C}$ be the $\mathbb{Q}$-isogeny class of $E^d$, where $E$ is an elliptic curve of conductor $p$ with some rational $\ell$-torsion. Then $\mathbb{Q}(A[\ell]) = \mathbb{Q}(\boldsymbol{\mu}_{2\ell})$ for all abelian varieties $A$ in $\mathcal{V}_\ell(\mathcal{C})$.*

*Proof.* We may assume that $E = E_\ddagger$ is $\ell$-maximal. Put
$$\mathcal{M}_2 = \mathcal{M}_2(E^d, v) = \oplus_1^d \mathcal{M}_2(E, v).$$
In view of (6.3), we have $\sigma\tau(z) = \tau(z) + s\,z$ for all $z \in \mathcal{M}_2$.

Suppose $A \in \mathcal{C}$ and let $\varphi : E^d \to A$ be the corresponding isogeny. The basic assumptions (C1), (C2) of section 4 clearly are satisfied by $E^d$. Condition (C3) holds because $\rho_{E^d} = \oplus \rho_E$ is a sum of $d$ copies of a fixed representation of $\mathrm{GL}_2$-type. Since these conditions are isogeny invariant, they also hold for $A$. Suppose, in addition, that $A$ is $\ell$-maximal. Then Proposition 4.12 implies that $\mathbb{T}_\ell(A) = \mathcal{M}'_2 \oplus \tau(\mathcal{M}'_2)$, where $\mathcal{M}'_2 = \mathcal{M}_2(A, v) = \mathcal{M}_1(A, v)$.

By construction, $\sigma$ acts trivially on $\mathcal{M}'_2$. To examine the action of $\sigma$ on $\tau(\mathcal{M}'_2)$, recall from the start of section 4 that $\mathcal{M}'_2$ contains $\varphi(\mathcal{M}_2)$ with finite index, say $n$. Given $z' \in \mathcal{M}'_2$, we may therefore find $z \in \mathcal{M}_2$ such that $n\,z' = \varphi(z)$. It follows that
$$n\sigma\tau(z') = \varphi(\sigma\tau(z)) = \varphi(\tau(z) + s\,z) = n\tau(z') + ns\,z'.$$

But $\mathbb{T}_\ell(A)$ is torsion-free, so $\sigma\tau(z') = \tau(z') + s\,z'$. Since $s \equiv 0 \bmod \ell$, we find that $\sigma$ acts trivially on the first layer $\tau(\bar{\mathcal{M}}_2')$ of $\tau(\mathcal{M}_2')$ and therefore on all of $A[\ell] = \bar{\mathcal{M}}_2' \oplus \tau(\bar{\mathcal{M}}_2')$. Hence $\mathbb{Q}(A[\ell]) = \mathbb{Q}(\boldsymbol{\mu}_{2\ell})$ by Lemma 3.6 and Corollary 4.5. □

**Lemma 6.8.** *Suppose that $\mathcal{C}$ is a $\mathbb{Q}$-isogeny class of abelian varieties of dimension $d$ satisfying (C1), (C2), (C3) and that $\mathbb{Q}(A[\ell]) \subseteq \mathbb{Q}(\boldsymbol{\mu}_{2\ell})$ for each $A$ in $V_\ell(\mathcal{C})$. Then there exists a $\mathbb{Z}_\ell$-submodule $X$ of $\mathbb{T}_\ell(A)$ of rank 2 and stabilized by $G_\infty$ such that $\mathbb{T}_\ell(A) \simeq X^d$ as $G_\infty$-modules.*

*Proof.* Theorem 5.4 implies that $A$ has totally toroidal reduction modulo $p$, so $\mathcal{M}_2 = \mathcal{M}_2(A,v) = \mathcal{M}_1(A,v)$. For all $A \in \mathcal{V}_\ell(\mathcal{C})$, we have $\bar{\mathcal{M}}_2 \cap \tau(\bar{\mathcal{M}}_2) = 0$ by (4.13) and $\mathbb{T}_\ell(A) = Y_A = \mathcal{M}_2 \oplus \tau(\mathcal{M}_2)$ by Proposition 4.12.

Choose a set of free generators $m_1, \ldots, m_d$ for $\mathcal{M}_2$. For $j = 1, \ldots, d$, define $X^{(j)}$ to be $\mathbb{Z}_\ell$-span of $m_j$ and $\tau(m_j)$. Each $X^{(j)}$ is a pure submodule of $\mathbb{T}_\ell(A)$ of rank 2 and $\mathbb{T}_\ell(A) = \oplus_{j=1}^d X^{(j)}$.

We wish to use Lemma 4.6 to show that $X = X^{(j)}$ is a $G_\infty$-module. Take $\mathcal{H}$ to be the group generated by $\tau$, so $X$ certainly is an $\mathcal{H}$-module by (C3). Suppose $X_n = X + \ell^n \mathbb{T}_\ell(A)$ is a $G_\infty$-module and let $\varphi : A \to A'$ be the $\mathbb{Q}$-isogeny whose kernel is $\kappa = (X + \ell^n \mathbb{T}_\ell(A))/\ell^n \mathbb{T}_\ell(A) \simeq \mathbb{Z}/\ell^n \oplus \mathbb{Z}/\ell^n$.

Write $\bar{\mathcal{M}}_2^{(n)}$ for the projection of $\mathcal{M}_2$ to $\mathbb{T}_\ell(A)/\ell^n \mathbb{T}_\ell(A) \simeq A[\ell^n]$. It is clear that $\kappa \cap \bar{\mathcal{M}}_2^{(n)}$ is isomorphic to one copy of $\mathbb{Z}/\ell^n$, generated by the coset of $m_j$. Hence $|\kappa \cap \bar{\mathcal{M}}_2^{(n)}| = \ell^n$ and we have $|\Phi_{A'}(\bar{\mathbb{F}}_p)|_\ell = |\Phi_A(\bar{\mathbb{F}}_p)|_\ell$ by Lemma 4.3. Therefore $A'$ also is $\ell$-maximal. Then by assumption, $\mathbb{Q}(A'[\ell]) \subseteq F = \mathbb{Q}(\boldsymbol{\mu}_{2\ell})$. Hence the subgroup $\mathcal{N}$ defined in Lemma 4.6 contains $\mathrm{Gal}(L_\infty/F)$. But $\mathrm{Gal}(L_\infty/F)$ and $\mathcal{H}$ certainly generate $G_\infty$. We may conclude from Lemma 4.6 that $X$ is a $G_\infty$-module.

The following standard argument now shows that the $X^{(j)}$'s are isomorphic as $G_\infty$-modules. Reasoning as above, we find that for each $j \neq 1$, the $\mathbb{Z}_\ell$-submodule of $\mathbb{T}_\ell(A)$ of rank 2 spanned by $m_1 + m_j$ and $\tau(m_1 + m_j)$ also is stabilized by $G_\infty$. But then the matrix representation of $G_\infty$ afforded by the $\mathbb{Z}_\ell$-span of $m_1$ and $\tau(m_1)$ must be identical to the representation afforded by the $\mathbb{Z}_\ell$-span of $m_j$ and $\tau(m_j)$. It now follows that there is a $G_\infty$-isomorphism $\mathbb{T}_\ell(A) \simeq X^d$, where $X$ is the $\mathbb{Z}_\ell$-module spanned by $m$ and $\tau(m)$ for any choice of $m \in \mathcal{M}_2$ such that $m \notin \ell \mathcal{M}_2$. □

For the proof of our final result, we need to observe that $\ell$-maximality is preserved for certain products of abelian varieties. Note that $A = B \times C$ satifies (C1), (C2), (C3) of section 4 if and only if each factor does.

**Lemma 6.9.** *Suppose $A$ satisfies (C1), (C2), (C3). If $A = B \times C$ is the product of $\ell$-maximal abelian varieties, then $A$ also is $\ell$-maximal.*

*Proof.* To show $A$ is $\ell$-maximal, it suffices to consider $\mathbb{Q}$-isogenies $\varphi : A \to A'$ of $\ell$-power degree. Set $\kappa = \mathrm{Ker}\,\varphi$ and assume $\kappa$ has exponent $\ell^n$. Theorem 5.4 implies that $A$ has totally toroidal reduction at $p$, so $\mathcal{M}_1(A,v) = \mathcal{M}_2(A,v)$ and similarly for $B$ and $C$. Moreover, inside $\mathbb{T}_\ell(A) = \mathbb{T}_\ell(B) \times \mathbb{T}_\ell(C)$, we have $\mathcal{M}_1(A,v) = \mathcal{M}_1(B,v) \times \mathcal{M}_1(C,v)$. Write $\bar{\mathcal{M}}_2^{(n)}$ for projection to $\mathbb{T}_\ell/\ell^n \mathbb{T}_\ell$ and let

$$\kappa_0 = \kappa \cap \bar{\mathcal{M}}_2^{(n)}(A,v) = \{\,(\bar{b},\bar{c}) \in \kappa \mid b \in \mathcal{M}_2(B,v),\ c \in \mathcal{M}_2(C,v)\,\} \subset A[\ell^n].$$

In view of (4.13), we have

$$\bar{\mathcal{M}}_2^{(n)}(B,v) \cap \tau(\bar{\mathcal{M}}_2^{(n)}(B,v)) = \bar{\mathcal{M}}_2^{(n)}(C,v) \cap \tau(\bar{\mathcal{M}}_2^{(n)}(C,v)) = 0,$$

whence $\kappa_0 \cap \tau(\kappa_0) = 0$. But $\kappa$ is a Galois module, so it contains $\kappa_0 + \tau(\kappa_0)$ and therefore $|\kappa| \geq |\kappa_0| |\tau(\kappa_0)| = |\kappa_0|^2$. It follows from Lemma 4.3 that $|\Phi_A(\bar{\mathbb{F}}_p)|_\ell \geq |\Phi_{A'}(\bar{\mathbb{F}}_p)|_\ell$. Hence $A$ is $\ell$-maximal. □

**Theorem 6.10.** *Let $\mathcal{C}$ be a $\mathbb{Q}$-isogeny class of semistable abelian varieties of dimension d with bad reduction only at p. Assume $\ell$ is a regular prime and $\mathbb{Q}(A[\ell]) \subseteq \mathbb{Q}(\boldsymbol{\mu}_{2\ell})$ for every $A$ in $\mathcal{V}_\ell(\mathcal{C})$. Suppose the action of $\tau$ on $\mathbb{T}_\ell(A)$ satisfies a polynomial of degree at most 2. Then $\ell = 2, 3$ or 5 and there exists an elliptic curve $E$ of conductor $p$ with some rational $\ell$-torsion such that $\mathcal{C}$ is the $\mathbb{Q}$-isogeny class of $E^d$.*

*Proof.* We may apply Lemma 6.8 to $A \in \mathcal{V}_\ell(\mathcal{C})$ to obtain a $\mathbb{Z}_\ell$-submodule $X$ of $\mathbb{T}_\ell(A)$ of rank 2 stabilized by $G_\infty$, such that $\mathbb{T}_\ell(A) \simeq X^d$ is a $G_\infty$-isomorphism. According to the Tate conjecture, proved by Faltings ([Sch]), we have

$$\mathrm{End}_\mathbb{Q}(A) \otimes \mathbb{Z}_\ell \simeq \mathrm{End}_{G_\infty} \mathbb{T}_\ell(A).$$

But the commutant $\mathrm{End}_{G_\infty} \mathbb{T}_\ell(A)$ clearly contains $M_d(\mathbb{Z}_\ell)$. It follows that the $\mathbb{Z}_\ell$-rank of $\mathrm{End}_\mathbb{Q}(A) \otimes \mathbb{Z}_\ell$ is at least $d^2$.

Suppose $A$ is simple, so that $D = \mathrm{End}_\mathbb{Q}(A) \otimes \mathbb{Q}$ is a division algebra. We have shown above that $[D : \mathbb{Q}] \geq d^2$. But $D$ acts on the space of invariant holomorphic differentials $\Omega^1(A)$, so

$$d = \dim_\mathbb{Q} \Omega^1(A) = [D : \mathbb{Q}] \dim_D \Omega^1(A) \geq d^2.$$

Hence $d = 1$ and $A$ is a Neumann-Setzer ($\ell = 2$) or Miyawaki ($\ell = 3$ or 5) elliptic curve.

If $A' \in \mathcal{V}_\ell(\mathcal{C})$ is not $\mathbb{Q}$-simple, we may find proper abelian subvarieties $B'$ and $C'$, such that $A'$ is $\mathbb{Q}$-isogenous to $B' \times C'$. Certainly $B'$ and $C'$ inherit properties (C1), (C2), (C3) from $A'$. Choose $B$ and $C$ to be $\ell$-maximal in the isogeny class of $B'$ and $C'$ respectively. Then $A'$ is isogenous to $A = B \times C$ and Lemma 6.9 shows that $A$ is in $\mathcal{V}_\ell(\mathcal{C})$. By assumption, the $\ell$-division field of $A$ therefore is contained in $\mathbb{Q}(\boldsymbol{\mu}_{2\ell})$ and so *a fortiori* the same is true of the $\ell$-division fields of $B$ and $C$. Thus we may complete the proof by induction. □

**Corollary 6.11.** *Let $\mathcal{C}$ be a $\mathbb{Q}$-isogeny class of semistable abelian varieties of dimension d with bad reduction only at p and suppose $\mathrm{Gal}(\mathbb{Q}(A[2])/\mathbb{Q})$ is abelian for all $A \in \mathcal{V}_2(\mathcal{C})$. Then $\mathcal{C}$ is the isogeny class of $E^d$, where $E$ is a Neumann-Setzer curve of conductor $p$.*

*Proof.* Theorem 1.1 tells us that $A$ is totally toroidal, so for $A \in \mathcal{V}_\ell(\mathcal{C})$, we have $A[2] = \bar{\mathcal{M}}_2 \oplus \tau(\bar{\mathcal{M}}_2)$ by (4.13) and equality of dimension. Since $\sigma$ acts trivially on $\mathcal{M}_2$ and the restrictions of $\sigma$ and $\tau$ to $\mathbb{Q}(A[2])$ commute by assumption, we find that $\sigma$ acts trivially on $A[2]$. It follows that $p$ is unramified in $\mathbb{Q}(A[2])$, whence $\mathbb{Q}(A[2]) = \mathbb{Q}(i)$ by Corollary 4.5. Our claim is now implied by Theorem 6.10. □

The following two examples show that the $\ell$-division field of a *simple* semistable abelian variety of dimension greater than 1 with one place of bad reduction may be an $\ell$-group. Compare with Theorem 6.10 and Theorem 1.4.

From a model for the modular curve $X_0(41)$ given by Weber [Web], one finds the minimal model

$$y^2 + (x^4 + 2x^3 - x)y = -3x^6 - 16x^5 - 29x^4 - 14x^3 + 13x^2 + 9x - 4.$$

Its Jacobian $J_0(41)$ has real multiplications by the maximal order in the totally real cubic field of discriminant 148 and bad reduction only at 41, where it is totally toroidal. The 2-division field of $J_0(41)$, can be seen to be cyclic over $K = \mathbb{Q}(\sqrt{-41})$ and dihedral over $\mathbb{Q}$, as we expect from Proposition 3.7. Compare with Proposition 6.6(a). Since the class number of $K$ is 8, the maximal $(2, 41)$-controlled extension has degree 32 over $\mathbb{Q}$. However $[L_0 : \mathbb{Q}] = 16$.

A minimal model for $X_0(31)$ is given by

$$y^2 + (x^3 + x^2 + 2x + 1)y = -x^4 - 2x^3 - 3x^2 - 2x - 1.$$

Its Jacobian variety $A = J_0(31)$ admits real multiplications by the ring of integers of $\mathbb{Q}(\sqrt{5})$ and is the full Eisenstein quotient at 5, cf. [Ma]. For the endomorphism $\pi = \sqrt{5}$, we have $A[\pi] = \mathbb{Z}/5\mathbb{Z} \oplus \boldsymbol{\mu}_5$ and so $H = \mathrm{Gal}(\mathbb{Q}(A[5])/\mathbb{Q}(\boldsymbol{\mu}_5))$ is an elementary abelian 5-group. One can check that its $\mathbb{F}_5$-rank is 3, in contrast to Proposition 6.6(c). Using the Remark after Lemma 3.4, we find that $\mathbb{Q}(A[5])$ is the maximal $(5, 31)$-controlled extension abelian over $\mathbb{Q}(\boldsymbol{\mu}_5)$.

(A. Brumer) DEPARTMENT OF MATHEMATICS, FORDHAM UNIVERSITY, BRONX, NY 10458
*E-mail address*: `brumer@fordham.edu`

(K. Kramer) DEPARTMENT OF MATHEMATICS, QUEENS COLLEGE (CUNY), FLUSHING, NY 11367
*E-mail address*: `kramer@forbin.qc.edu`